\documentclass{article}

\usepackage{latexsym}
\usepackage{amssymb,amsfonts,amsmath}
\usepackage{stmaryrd}

\newcommand{\al}{\alpha}
\newcommand{\be}{\beta}

\newcommand{\eps}{\varepsilon}

\newcommand{\Ab}{{\bf A}}
\newcommand{\Cb}{{\bf C}}
\newcommand{\Wb}{{\bf W}}
\newcommand{\Ub}{{\bf U}}
\newcommand{\Tb}{{\bf T}}
\newcommand{\Sb}{{\bf S}}

\newcommand{\R}{{\bf R}}
\newcommand{\Z}{{\bf Z}}
\newcommand{\C}{{\bf C}}
\newcommand{\T}{{\bf T}}

\newcommand{\Q}{{\bf Q}}
\newcommand{\N}{{\bf N}}
\renewcommand{\H}{{\cal H}}

\newcommand{\Ac}{{\cal A}}
\newcommand{\Bc}{{\cal B}}
\newcommand{\Sc}{{\cal S}}
\newcommand{\Fc}{{\cal F}}
\newcommand{\Gc}{{\cal G}}
\newcommand{\ab}{{\bf a}}
\newcommand{\bb}{{\bf b}}
\newcommand{\cb}{{\bf c}}
\newcommand{\db}{{\bf d}}
\newcommand{\xb}{{\bf x}}

\newcommand{\Tau}{\mathcal{T}}

\newcommand{\norm}[1]{\|#1\|}
\newcommand{\ip}[2]{\langle#1,#2\rangle}

\newtheorem{Theorem}{Theorem}[section]
\newtheorem{Lemma}[Theorem]{Lemma}
\newtheorem{Corollary}[Theorem]{Corollary}

\newtheorem{Example}[Theorem]{Example}
\newtheorem{Remark}[Theorem]{Remark}

\title{A Noncommutative Wiener Lemma and A Faithful Tracial State
 on Banach Algebras of Time-Frequency Shift Operators}
\author{Radu Balan \\ Siemens Corporate Research \\ 
755 College Road East \\ Princeton, NJ 08540 \\ e-mail: {\em
radu.balan@siemens.com} }

\date{\today}

\begin{document}

\maketitle

\begin{abstract}
In this paper we analyze the Banach *-algebra of time-frequency shifts
with absolutely summable coefficients. We prove a noncommutative version
of the Wiener lemma. We also construct a faithful tracial state on this 
algebra which implies the algebra contains no compact operators. 
As a corollary we obtain
a special case of the Heil-Ramanathan-Topiwala conjecture regarding
linear independence of finitely many time-frequency shifts of one
$L^2$ function.
\end{abstract}

\section{Introduction\label{sec1}}

The Time-Frequency representation of the Heisenberg group has received a
lot of attention for the past 20 years with the advent of Gabor
analysis. Many methods and techniques have been developed and a rich
body of results has been obtained. For a nice account of such results
we refer the reader to the excellent monography \cite{charlybook}.

For $t\in\R^d$ and $\omega\in\R^d$ we denote by $S_t$ the time shift operator,
 by $M_\omega$ the frequency shift operator, and by $U_{t,\omega}$  the 
time-frequency shift operator defined, respectively, by:
\begin{eqnarray}
S_t:L^p(\R^d)\rightarrow L^p(\R^d) & , & S_tf(x)=f(x-t) \\
M_\omega:L^p(\R^d)\rightarrow L^p(\R^d) & , & M_\omega f(x)=e^{i\omega x}f(x) 
\\
U_{t,\omega}:L^p(\R^d)\rightarrow L^p(\R^d) & , & U_{t,\omega}f(x)=
M_\omega S_t f(x) = e^{i\omega x}f(x-t)
\end{eqnarray}
For a Banach space $X$, we let $B(X)$ denote the Banach space of
bounded operators on $X$ with the usual operator norm.  Consider now
the Banach algebra (see \cite{Neimark,Rickart} for definition and
properties of Banach algebras) of time-frequency shifts on $L^2(\R^d)$ with
absolutely summable coefficients:
\begin{equation}\label{eq:1.1}
\Ac = \{T\in B(L^2(\R^d))~~|~~ T=\sum_{\lambda\in\R^{2d}}c_\lambda U_\lambda~,~
\norm{T}_{\Ac}:=\sum_{\lambda\in\R^{2d}}|c_\lambda|<\infty\}
\end{equation}
Note the following: (i) for any $T\in\Ac$, the support $supp(T)$ of its
generating sequence $c$, $supp(T)=\{\lambda\in\R^{2d}~,~c_\lambda\neq 0\}$ 
is always a countable set;
 (ii) this support is not assumed to have any lattice structure.


Wiener's lemma states that is a periodic function $f$ has an
absolutely convergent Fourier series and never vanishes, then its
reciprocal $1/f$ has also an absolutely convergent Fourier series. In
Banach algebras language the same result can be restated as follows.
Consider the Banach algebra
\begin{equation}
 {\Ac}_W = \{T\in B(L^2(\R^d))~~|~~T =\sum_{n\in\Z^d}c_n S_n~,~
\norm{T}_{{\Ac}_W}:=\sum_n |c_n| 
\}
\end{equation}
which is a $*$-subalgebra of $B(L^2(\R^d))$. Then Wiener's lemma asserts
that, if $T\in\Ac_W$ is invertible in $B(L^2(\R^d))$, then
$T^{-1}\in{\Ac}_W$. In Naimark's terminology (\cite{Neimark}),
$(\Ac_W,B(L^2(\R^d)))$ is said a {\em Wiener pair}, in
 \cite{baskakov97} 
 $\Ac_W$ is called a {\em full algebra} in $B(L^2(\R^d))$, whereas
\cite{Barnes00} calls such algebras {\em inverse closed}.

Many generalizations of this result appeared in literature. We will
mention here four extensions that set the context of our 
results. 

In the context of almost periodic functions, one obtains the following
result. Consider the Banach algebra
\begin{equation}
 \Ac_{AP}=\{T\in B(L^2(\R^d))~~|~~ T=\sum_{t\in\R^d} c_t S_t~~|~~
\norm{T}_{\Ac_{AP}}:=\sum_{t\in\R^d}|c_t|<\infty \}
\end{equation}
of absolutely summable linear combinations of arbitrary real
shifts. Then using the Bohr compactification of $\R^d$ and the
Gelfand's method of proving Wiener's lemma (\cite{Loomis}), it follows that
$\Ac_{AP}$ is an inverse closed algebra in $B(L^2(\R^d))$. More explicitely, if
$T=\sum_{\lambda\in\Lambda}c_{\lambda}S_{\lambda}$ for some countable
subset $\Lambda\subset\R^d$, and
$\sum_{\lambda\in\Lambda}|c_{\lambda}|<\infty$, then
$T^{-1}=\sum_{\sigma\in\Sigma}d_{\sigma}S_{\sigma}$, for some (in
general) other countable subset $\Sigma\subset\R^d$, and absolutely
summable complex coefficients
$\sum_{\sigma\in\Sigma}|d_{\sigma}|<\infty$. See also \cite{rodman98} for
an extension to matrix valued almost periodic functions.

In the context of time-frequency analysis Gr\"{o}chenig and Leinert
\cite{grochJAMS} obtained a discrete noncommutative Wiener lemma as follows.
 Fix $\al,\be>0$ and a subexponential weight $v$
(see section \ref{sec2}).  We let $\Ac_{v;GL}(\al,\be)$ denote the operator
algebra
\begin{equation}
\label{eq1.3}
\Ac_{v;GL}(\al,\be) = \{ T\in B(L^2(\R^d))~~|~~T=\sum_{k,l\in\Z^d}a_{k,l} 
S_{\alpha k}M_{\beta l}~~|~~\norm{T}_{\Ac_v}:=\norm{\ab}_{1,v}:=
\sum_{k,l\in\Z^d}v(k,l)|a_{k,l}| \}
\end{equation}
where $\norm{\ab}_{1,v}$ is the $v$-weighted $l^1$ norm of $\ab$.  
Then, using Ludwig's theorem on
symmetric group algebras of nilpotent groups (\cite{Ludwig}), the
authors proved in \cite{grochJAMS} (Theorem 3.1) that
$\Ac_{v:GL}(\al,\be)$ is inverse closed in $B(L^2(\R^d))$. The algebra
$\Ac_{v;GL}$ is naturally associated to the twisted convolution
algebra $(l^1_v(\Z^{2d}),\sharp_{\theta})$, where
\begin{equation}
(\ab \sharp_{\theta} \bb )(m,n)=\sum_{k,l\in\Z^d}a_{kl}b_{m-k,n-l}
e^{2\pi i\theta(m-k)\cdot l}
\end{equation}
In this setting, the above result is equivalent to saying (Theorem 2.14 in
\cite{grochJAMS}), if $\ab\in l^1_v(\Z^{2d})$ is so that the
convolution operator $L_{\ab}:l^2(\Z^{2d})\rightarrow l^2(\Z^{2d})$,
$\xb\mapsto L_{\ab}(\xb)=\ab\sharp_{\theta}\xb$, is invertible in
$B(l^2(\Z^{2d}))$, then its inverse is of the form $L_{\ab}^{-1}=L_{\bb}$
for some $\bb\in l^1_v(\Z^{2d})$. The proof of their results rely heavily on
fairly abstract results on group algebras of locally compact
nilpotent groups (\cite{Ludwig},\cite{Hulanicki72}). As the authors point
out, analyzing spectral properties of group algeabras is not usually an easy
business. 

Again in the context of time-frequency analysis, Gr\"{o}chenig in
\cite{groch05} translated a result by J.Sj\"{o}strand \cite{sjostrand95} using
the modulation space $M^{\infty,1}_v$. The operator algebra 
\begin{equation}
\label{eq1.3a}
\Sc_{v}=\{T\in B(L^2(\R^d))~~|~~T=\int_{\R^{2d}}d\lambda\sigma(\lambda)
U_{\lambda}~~,~~\norm{T}_{\Sc_v}:=\norm{\sigma}_{M^{\infty,1}_v}:=
\int_{\R^d}dq\,sup_{p\in\R^d}|V_{\gamma}\sigma(q,p)|v(q,p)\}
\end{equation}
is shown to be inverse closed in $B(L^2(\R^d))$.

Another example of a full algebra is furnished by the Baskakov
class of matrices that have some off-diagonal decay. 
In \cite{baskakov97} Baskakov proves the Banach algebra:
\begin{equation}
B_{v}=\{ A=(A_{m,n})_{m,n\in\Z^d}\in B(l^2(\Z^d))~~|~~\norm{A}_{B_v}:=
\sum_{k\in\Z^d}v(k)\sup_{m\in\Z^d}|A_{m,m-k}| <\infty\}
\end{equation}
with $v$ a subexponential weight is inverse closed and also obtains
estimates of the entries of the inverse matrix. The unweighted version
of this result had been proved in \cite{gohberg89} as reported to me
by T.Strohmer \cite{thomas05}. These results have been obtained also
independently by Gr\"{o}chening and Leinert in
\cite{GrochLeinert05} using a Banach algebra technique. In
\cite{excess3} we used Baskakov's result to establish localization 
results for Gabor like frames.

In this paper we extend previously known results to the Banach
algebra (\ref{eq:1.1}) and its weighted version (\ref{eq:Av}). Beside the
intrinsec interest of a new Wiener type lemma, we are motivated
by two problems of time-frequency analysis. One problem is
 {\em the Heil-Ramanathan-Topiwala (HRT) conjecture}, the other
problem relates to the time-frequency analysis of communication
channels.  

A second contribution of this paper is the explicit construction of a 
faithful tracial
state on $\Ac$ that yields several consequences.  In
particular we show that $\Ac$ does not contain any compact operator,
from where we obtain as a corrolary a partial answer to the HRT
conjecture. On the other hand we prove extensions of the Paley-Wiener
theorem for this algebra.

The HRT conjecture (see \cite{heilHRT}) states that finitely many distinct
time-frequency shifts of one $L^2$ function, are linearly independent
(over $\C$). This means, for any finite subset of $\R^{2d}$,
$\Lambda\subset\R^{2d}$, and $g\in L^2(\R^d)$,
\begin{equation}
\sum_{\lambda\in\Lambda}c_{\lambda}U_\lambda g=0 \Rightarrow c_{\lambda}=0,
\forall \lambda\in\Lambda 
\end{equation}
When $\Lambda$ is a subset of a lattice, the claim was positively
proved by Linnell in \cite{linnel}, however the general case, as far
as we know, is still open. The problem can be recast into a spectral
analysis problem.  More specifically the HRT conjecture is equivalent
to proving that for any finite subset $\Lambda\subset\R^{2d}$, and
complex numbers $(c_\lambda)_{\lambda\in\Lambda}$, the bounded
operator $T=\sum_{\lambda\in\Lambda}c_{\lambda}U_{\lambda}$ has no
pure point spectrum. Motivated by this problem, one is naturally led
to an algebra of type (\ref{eq:1.1}). Our result in this paper (Theorem
 \ref{t4}) is one step toward analyzing spectral properties of such operators.
 Here we prove that any $T$ of this form cannot have 
isolated eigenvalues of finite multiplicity.

In communication theory, a multipath time-varying communication
channel is modeled as a linear superposition of time-frequency shifts
(see \cite{thomasS}). Often the channel model contains finitely many
time-frequency shifts, or infinitely many but fast-decaying
coefficients, so naturally, the channel transfer operator is in
algebra $\Ac_v$.  One problem is channel equalization (or
deconvolution) by which one has to invert the channel transfer
operator. Assuming this operator is invertible on the space of finite
energy signals, then our result says the inverse is also a
superposition of time-frequency-shifts, with absolutely summable
coefficients.  The coefficients decaying rate gives the convergence
rate of finite approximation methods. In this context our results
(Theorems \ref{t2.1}, \ref{t5}) give estimates of this decay. We also
obtrain necessary and sufficient conditions for operators in $\Ac$
to have bounded support (Theorem \ref{t6}).

The organization of this paper is the following. In Section \ref{sec2}
we present algebra constructions; In Section
\ref{sec2b} we state our main results. In Section \ref{sec3} we
connect our approach to prior literature, and in Section \ref{sec4} we
prove these results.

Throughout this paper we use the following notations: For a set $I$,
$|I|$ denotes the cardinal of set $I$ (i.e. the number of points that $I$
contains); for $x\in\R^d$, $|x|$ denotes its $\infty$-norm, and
$\norm{x}$ the Euclidian ($l^2$) norm; $B_r(x)$ denotes the closed
ball of radius $r$ centered at $x$ with respect to norm $|\cdot|$,
 $B_r(x)=\{y\in\R^d~|~|x-y|\leq r\}$; thus
$\{B_1(n)~;~n\in\Z^d\}$ forms a covering (but not disjoint) partition of 
$\R^d$; $E_r(x)$ denotes the Euclidian closed ball of radius
$r$ centered at $x$, $E_r(x)=\{y\in\R^d~|~\norm{x-y}\leq r\}$; $\Fc$ denotes
the unitary Fourier transform with the following normalization:
\[ f \mapsto \Fc f(\omega) =  
\frac{1}{(2\pi)^{d/2}}\int_{\R^d} e^{-i\omega x}f(x)\,dx \]
We will frequently use $\lambda,\mu$ to denote time-frequency points in 
$\R^{2d}$, e.g. $\lambda=(t,\omega)$ of components $t,\omega\in\R^d$.

\section{Twisted convolution, Heisenberg group representation, and Banach 
algebras\label{sec2}}


\subsection{Banach Algebras}
A {\em Banach algebra } $A$ is a normed algebra which is closed with
respect to its norm $\norm{\cdot}_A$ and satisfies
$\norm{xy}_A\leq\norm{x}_A\norm{y}_A$ for all $x,y\in A$. Throughout this
paper all algebras have unit. In general, when an algebra
does not have an identity element, one formally adds a unit, and defines the
inverse with respect to this element. The {\em
resolvent} of an element $x\in A$ is the set of complex numbers
$\lambda\in\C$ so that $\lambda 1-x$ is invertible in $A$. Its
complement is called {\em spectrum} and is denoted by
$sp_A(x)$. It is always a closed set included into the ball or radius
$\norm{x}_A$ centered at the origin. The largest absolute value of
elements of the spectrum is called the {\em spectral radius}, denoted
$r_A(x)$, $r_A(x):=max_{\lambda\in sp_A(x)}|\lambda|$. Note
$r_A(x)\leq \norm{x}_A$. The spectral radius can be computed using
Gelfand's formula:
\begin{equation}\label{eq:Gelfand}
r_A(x)=\lim_{n\rightarrow\infty}(\norm{x^n}_A)^{1/n}
\end{equation}
A map $*$ of $A$ is called {\em involution}, $x\mapsto x^*$, if
$(ax+by)^*=\bar{a}x^*+\bar{b}y^*$, $(x^*)^*=x$, and
$\norm{x^*}_A=\norm{x}_A$ for all $a,b\in\C$ and $x,y\in A$. A Banach
algebra with an involution is called a {\em Banach $*$-algebra}.

If the norm $\norm{\cdot}_A$ satisfies $\norm{x^*x}_A=\norm{x}_A^2$, 
then $A$ is called a {\em $C^*$ algebra}.

Typically on a $*$-algebra $A$ we will have two norms:
$\norm{\cdot}_A$ with respect to which $A$ is closed, hence a
Banach $*$-algebra, and $\norm{\cdot}_B$ that satisfies
$\norm{x^*x}_B=\norm{x}_B^2$ and with respect to which $A$ is not
closed. The completion $B$ of $A$ with respect to this latter norm
$\norm{\cdot}_B$ is a $C^*$ algebra. Note $A\subset B$, 
hence $\norm{x}_B\leq\norm{x}_A$, for all $x\in A$. 

We call $A$ an {\em inverse closed} algebra in $B$, if any element 
$x\in A$ that is invertible in $B$, is invertible in $A$, $x^{-1}\in A$. 
Neimark in \cite{Neimark} calls $(A,B)$ a {\em Wiener
pair}, whereas Baskakov in \cite{baskakov97} calls $A$ a {\em full} 
algebra in $B$.

\subsection{Weighted Algebras} 
 
In this paper a weight $v$ is a nonnegative and radially non-decreasing 
function on $\R^d$ so that $v(0)=1$ and $v(-x)=v(x)$. 
Let $w:\R^{+}\rightarrow\R^{+}$ be the function $w(r)=max_{\norm{x}=r}v(x)$.
We define the following (see also \cite{GrochLeinert05}):

(a) The weight $v$ is said {\em submultiplicative} if satisfies
\begin{equation}
\label{eq:wsub}
v(x+y)\leq v(x)v(y)
\end{equation}

(b) The weight $v$ is said to satisfy the {\em GRS (Gelfand-Raikov-Shilov)
 condition} if
\begin{equation}
\label{eq:GRS}
\lim_{n\rightarrow\infty}(w(nr))^{1/n} = 1~~,~~\forall r\geq 0
\end{equation}

(c) The weight $v$ is called {\em admissible} if it is submultiplicative
and satisfies the GRS condition. 

\begin{Example}\cite{GrochLeinert05}

The typical examples of admissible weights are the polynomial weights,
 $v(x)=(1+\norm{x})^s$ for some $s\geq 0$, and the subexponential weights,
 $v(x)=e^{\alpha \norm{x}^\beta}$, for some $\alpha\geq 0$ and $0<\beta<1$.
More generally, the following is also an admissible weight 
(see \cite{GrochLeinert05}), $v(x)=e^{\alpha\norm{x}^\beta}(1+\norm{x})^s
log^t(e+\norm{x})$, where $\alpha,s,t\geq 0$, $0<\beta<1$.

Note the exponential weight $v(x)=e^{\alpha\norm{x}}$ with $\alpha>0$
 is not admissible. It is submultiplicative, but does not satisfy 
the GRS condition.
\end{Example}

Throughout this paper all the weights are assumed at least submultiplicative.
Except for Lemma \ref{l1} and the proof of Theorem \ref{t5}, all weights
considered in the rest of the paper paper are admissible. 
In Lemma \ref{l1} we consider the less
restrictive submultiplicative weights to cover the case of exponential
weights needed in the proof of Theorem \ref{t5}.

For a weight $v$ we denote by $l^1_v(\R^{n})$ (or just $l^1_v$ when no
danger of confusion) the space of functions $c:\R^n\rightarrow\C$ so that
$\norm{c}_{l^1_v}:=\sum_{x\in\R^n}v(x)|c(x)|<\infty$. For $0<p<\infty$
we let $l^p(\R^n)$
 (or merely $l^p$, when no danger of confusion) denote the space of functions
$c:\R^n\rightarrow\C$ so that $\norm{c}_p:=(\sum_{x\in\R^n}
|c(x)|^p)^{1/p}<\infty$. We will frequently use the notation $c_x=c(x)$.  
The support of $c$ is defined by $supp(c)=
\{\lambda\in\R^n~|~c_\lambda\neq 0\}$, and for any $c$ in $l^p$ with 
$0<p<\infty$  or $l^1_v$ it is always a finite or countable subset 
of $\R^n$. For $p=\infty$, $l^\infty$ represents the set of bounded functions
on $\R^n$, not necessarily of finite or countable support, and
$\norm{c}_{\infty}=\sup_x |c_x|$. For
any $p\geq 1$, $l^p$ is a Banach space with $\norm{\cdot}_p$ norm.
Note $l^p(\R^n)$ is not separable for any $p$. In particular 
$l^2(\R^n)$ does not have a countable orthonormal basis. $l^p(\R^n)$ and
$l^1_n(\R^n)$ are the corresponding $L^p$ and $L^1_v$ spaces for
 $\R^n$ endowed with discrete topology.

We denote by $\Ac_v$ the algebra of time-frequency operators 
\begin{equation}
\label{eq:Av}
\Ac_v=\{ T=\sum_{\lambda}c_\lambda U_\lambda ~~;~~\norm{T}_{\Ac_v}:=
\sum_\lambda v(\lambda)|c_\lambda|<\infty \}
\end{equation}
that is the algebra of time-frequency shifts whose coefficients are
$l^1_v$ summable.
Clearly this is a subalgebra of the bounded operators $B(L^2(\R^d))$.
Furthermore, for every $T\in\Ac_v$, $\norm{T}_{B(L^2(\R^d))}\leq
\norm{T}_{\Ac_v}$. Thus the $C^*$ algebra $\Ab$ obtained by closing
any one of $\Ac_v$ with respect to the operator norm 
$\norm{\cdot}_{B(L^2(\R^d))}$ includes $\Ac$ and hence every $\Ac_v$.

\section{Main Results\label{sec2b}}

Previous section set the framework for stating our results. In the
following $v$ denotes a subexponential weight. In particular $v$ can
be the constant function $v=1$ (the unweighted case).

\begin{Theorem}[Spectral Invariance]\label{t2}
Assume $T=\sum_{\lambda}c_\lambda U_\lambda \in\Ac_v$. Then
the spectral radia with respect to the two algebras $B(L^2(\R^d))$ and
 $\Ac_v$ satisfy
\begin{equation}
\label{eq3.3.1}
r_{B(L^2(\R^d))}(T)=r_{\Ac_v}(T)
\end{equation}
\end{Theorem}

\begin{Theorem}[Wiener Lemma for TF Operators]\label{t1}
The algebra $\Ac_v$ is inverse closed in $B(L^2(\R^d))$. Explicitely
this means, if $T=\sum_{\lambda\in\Lambda}c_{\lambda}U_{\lambda}$ for
some $\cb\in l^1(\R^{2d})$ with $\Lambda=supp(\cb)$, and $T$ is
invertible in $B(L^2(\R^d))$, then there is $\db\in l^1(\R^{2d})$ with
$\Sigma=supp(\db)$ so that
$T^{-1}=\sum_{\sigma\in\Sigma}d_{\sigma}U_{\sigma}$.
\end{Theorem}

Immediate corollaries of this result are the following:
\begin{Corollary}\label{c0}
For any $T\in\Ac_v$ its spectrum with respect to the algebra $\Ac_v$ 
coincides to the spectrum with respect to the algebra $B(L^2(\R^d))$. 
Explicitely this means 
\begin{equation}
\label{eq3.3.2}
sp_{B(L^2(\R^d))}(T) = sp_{\Ac_v}(T)
\end{equation}
\end{Corollary}

\begin{Corollary}\label{c1}
Assume $T=\sum_{\lambda\in\R^{2d}}c_\lambda U_{\lambda}$ with $\sum_{\lambda}
|c_{\lambda}|<\infty$ is invertible in $B(L^2(\R^d))$. Then $T$ is
invertible in all $B(L^p(\R^d))$, with $0< p\leq\infty$.
\end{Corollary}

\begin{Corollary}\label{c2}
Let $T=\sum_{t\in\R^d}m_tS_t$ be a bounded invertible operator on $L^2(\R^d)$
so that $\sum_{t}\norm{m_t}_{AP}<\infty$. Then $T^{-1}=\sum_{t\in\R^d}
n_t S_t$ with $n_t\in AP$ so that $\sum_{t}\norm{n_t}_{AP}<\infty$.
\end{Corollary}

The following theorem gives an explicit estimate 
of the $\Ac_v$ norm of the inverse when
the invertible operator has finite support.
\begin{Theorem}[Norm of the Inverse]\label{t2.1}
Assume $T=\sum_{\lambda\in\Lambda}c_\lambda U_\lambda$ with 
$|\Lambda|=N<\infty$ and $R_0=\max_{\lambda\in\Lambda}\norm{\lambda}$. 
Assume $T$ is invertible in $B(L^2(\R^d))$, and
hence in $\Ac_v$ as well (by Theorem \ref{t1}). 
Denote $A=\norm{T^{-1}}^2_{B(L^2(\R^d))}$, $B=\norm{T}^2_{B(L^2(\R^d))}$, and
$\rho=\max(1,2R_0)$, and assume a polynomial weight $w(x)=C(1+x)^m$ for
some $C>0$ and $m\in\N$. Then
\begin{equation}\label{eq:Tninv}
\norm{T^{-1}}_{\Ac_v}\leq \frac{C\rho^m\norm{T}_{\Ac_v}}{A} (m+N)!
\left(\frac{A+B}{2A}\right)^{m+N}
\end{equation}
\end{Theorem}

Consider now $\Gc=\{g_{m,n;\al,\be}:=U_{\be n,2\pi\al m}g~|~m,n\in\Z^d\}$ 
a Gabor frame
for $L^2(\R^d)$, with $\al,\be>0$, $\al\be\leq 1$, and a dual Gabor frame
 (not necessarily the canonical dual frame)
 $\tilde{\Gc}=\{ \tilde{g}_{m,n;\al,\be}:=U_{\be n,2\pi\al m}\tilde{g}~|
~m,n\in\Z^d\}$. For details on Gabor frame theory we refer the reader 
to e.g. \cite{charlybook}. The following theorem
gives an explicit construction of the faithful tracial state:
\begin{equation}\label{eq:3.0}
\gamma:\Ac\rightarrow\C~~,~~\gamma(\sum_{\lambda}c_\lambda U_\lambda)=c_0
\end{equation}
This trace extends to $\Ab$, the completion of $\Ac$ with respect to the
 operator norm, which is a $C^*$ algebra.
\begin{Theorem}[Trace on $\Ab$]\label{t3}
For any $T\in \Ab$, 
\begin{equation}
\label{eq3.1}
\gamma(T) = \frac{1}{(\al\be)^d}\lim_{M,N\rightarrow\infty}\frac{1}{(2M+1)^d
(2N+1)^d}\sum_{|m|\leq M}\sum_{|n|\leq N}\ip{T g_{m,n;\al,\be}}{
\tilde{g}_{m,n;\al,\be}}
\end{equation}
is the faithful tracial state (\ref{eq:3.0}) on $\Ab$, independent 
of the choice of the Gabor frame $\Gc$.
\end{Theorem}

An immediate consequence of Theorem \ref{t4} is the following corollary
\begin{Corollary}\label{c3.1}
For any operator $T=\sum_\lambda c_\lambda U_\lambda\in\Ac$ 
\begin{equation}\label{eq:norms}
|c_\lambda|\leq \norm{c}_{\infty}\leq \norm{c}_2 \leq \norm{T}_{B(L^2(\R^d))}
\leq\norm{c}_1
\end{equation}
and
\begin{equation}\label{eq:coeffs}
c_\lambda = \gamma(U_\lambda^*T) = \lim_{M,N\rightarrow\infty}
\frac{1}{(\al\be)^d}\frac{1}{(2M+1)^d(2N+1)^d}\sum_{|m|\leq M}\sum_{|n|\leq N}
\ip{U_\lambda^*Tg_{m,n;\al,\be}}{\tilde{g}_{m,n;\al,\be}}
\end{equation}
\end{Corollary}

 As a corollary of this result we obtain that $\Ab$ (and hence $\Ac$ as well)
cannot contain compact operators:
\begin{Corollary}\label{c3}
Assume $T\in\Ab$ is a compact operator. Then $T=0$.
\end{Corollary}

Since for any finite set $\Lambda\in\R^{2d}$, and complex scalars
$(c_\lambda)$, the operator $T=\sum_{\lambda\in\Lambda}c_\lambda U_\lambda$
is in $\Ac$, we obtain the following theorem that gives a partial answer to
the HRT conjecture:
\begin{Theorem}\label{t4}
For any finite $\Lambda\subset\R^{2d}$ and complex scalars
$(c_\lambda)_{\lambda\in\Lambda}$, the operator
$T=\sum_{\lambda\in\Lambda} c_\lambda U_\lambda$ has no finite
multiplicity eigenvalue. Hence the pure point spectrum, if exists, can
only contain either eigenvalues with infinite multiplicity, or
eigenvalues that belong to the continuum part of the spectrum as well.
\end{Theorem}

For finite support operators as above, we can estimate the decay rate
of the coefficients of the inverse operator. In general such operators
are in $\Ac_v$ for any subexponential weight $v$. Hence the inverse
operator would have coefficients that are summable with such
weights. However, we can obtain more:

\begin{Theorem}\label{t5}
Let $T=\sum_{\lambda\in\Lambda}c_\lambda U_\lambda$ be an invertible
operator with $\Lambda\subset\R^{2d}$ a finite set. 
Then there is $\delta>0$ so that if $T^{-1}=\sum_{\mu\in\R^{2d}}
d_\mu U_\mu$ then
\begin{equation}\label{eq:3.3}
\sum_{\mu\in\R^{2d}} e^{\delta |\mu|}|d_\mu| < \infty
\end{equation}
\end{Theorem}
This result generalizes the classic statement (see
 e.g. \cite{zygmund}) that the reciprocal of a trigonometric
 polynomial that does not vanish on the unit circle has exponentially
 decaying Fourier coefficients. Our result is stronger than just
 simply saying the coefficients of the inverse operator should decay
 exponentially fast. It also controls the ``tail'' of the summation. To
 make this point more precise, we notice that (\ref{eq:3.3}) is equivalent to a
 Wiener amalgam type norm:
\begin{equation}
\sum_{n\in\Z^{2d}}e^{\delta |n|}\sum_{\mu\in B_1(n)}|d_\mu| < \infty
\end{equation}
where $B_1(n)$ is the ball of radius 1 centered at $n$,
$B_1(n)=\{\lambda\in\R^{2d}~|~|n-\lambda|\leq 1\}$.  
Hence there is a constant $C>0$ so that for all $R>0$,
\begin{equation}
\label{eq:3.3b}
\sum_{\mu\in\R^{2d},|\mu|\geq R}|d_\mu|<Ce^{-\delta R}
\end{equation}

Another equivalent statement of Theorem \ref{t5} is given by the
 Corollary \ref{c5.1} below. First we need to introduce a Banach space.
 Let us denote by $L^{2,\infty}(\R^d\times I)$ the mixed
norm Banach space
\begin{equation}\label{eq:L2inf}
L^{2,\infty}(\R^d\times I)
=\{f:\R^d\times I\rightarrow\C~;~\norm{f}_{2,\infty}^2:=
\sup_{y\in I}\int_{\R^d}|f(x,y)|^2dx<\infty \}
\end{equation}
where $I\subset\R^d$ is a compact neighborhood of the origin.
Then the unitary $U_\lambda$ extends from $L^2(\R^d)$ to $L^{2,\infty}(\R^d
\times I)$ simply by:
\begin{equation}
\Ub_\lambda f(x,y) = e^{i\ip{\omega}{x+iy}}f(x-t,y)~~,~~
\lambda=(t,\omega),~(x,y)\in\R^d\times I
\end{equation}
An operator $T=\sum_\lambda c_\lambda U_\lambda$ extends to
$ \Tb = \sum_\lambda c_\lambda \Ub_\lambda$ under some conditions.
Clearly all finite or compactly supported operators of $\Ac$ can be
extended to $L^{2,\infty}(\R^d\times I)$. Theorem \ref{t6} gives
necessary and sufficient conditions for such an extension to exist.

\begin{Corollary}\label{c5.1}
Let $T=\sum_{\lambda\Lambda}c_\lambda U_\lambda$ be a finitely supported
 invertible operator in $\Ac$. Then its inverse $T^{-1}=\sum_\mu d_\mu U_\mu$
extends to $L^{2,\infty}(\R^d\times I)$ to the inverse of the extension:
\begin{equation}
\label{eq:88}
\Tb^{-1} = \sum_\mu d_\mu \Ub_\mu
\end{equation}
for some compact neighborhood $I$ of the origin whose size depends on
the operator $T$.
\end{Corollary}

\begin{Theorem}\label{t6}
\begin{enumerate}
\item Assume $T=\sum_{\lambda\in\Lambda}c_\lambda U_\lambda\in\Ac$ so that
$\Lambda\subset\R^d\times E_\Omega(0)$ for some $\Omega>0$, where 
$E_\Omega(0)=\{\omega\in\R^d~|~\norm{\omega}\leq\Omega\}$. 
Then for any $\rho>0$, $T$ extends to 
$X_\rho=L^{2,\infty}(\R^d\times E_\rho(0))$ 
with operator norm bounded by:
\begin{equation}
\label{eq:3.6}
\norm{\Tb}_{B(X_\rho)}\leq C e^{\rho\Omega}
\end{equation}
for $C=\norm{T}_{\Ac}$.

\item Conversely, assume $T=\sum_{\lambda\in\Lambda}c_\lambda U_\lambda\in\Ac$
 can be extended to $X_\rho$ for all $\rho>0$, with a norm bounded as 
in (\ref{eq:3.6}),
 for some $C>0$ and $\Omega>0$ independent of $\rho$. 
Then $supp(T)\subset\R^d\times E_\Omega(0)$.

\item Assume $T=\sum_{\lambda\in\Lambda}c_\lambda U_\lambda\in\Ac$ so that
$\Lambda\subset E_D(0)\times\R^d$ for some $D>0$. Then for any $\rho>0$,
the operator $S=\Fc^* T\Fc$, where $\Fc$ denotes the Fourier transform,
 extends to $X_\rho$ with operator norm bounded by
\begin{equation}
\label{eq:3.7}
\norm{\Sb}_{B(X_\rho)}\leq C e^{\rho D}
\end{equation}
for $C=\norm{T}_{\Ac}$.

\item Conversely, assume the Fourier conjugate $\Fc^* T\Fc$ of
 $T=\sum_{\lambda\in\Lambda}c_\lambda U_\lambda\in\Ac$ can be extended
to $X_\rho$ for all $\rho>0$ with a norm bounded as in (\ref{eq:3.7})
with $C$ and $D$ independent of $\rho$. 
Then $supp(T)\subset E_D(0)\times\R^d$. 
\end{enumerate}
\end{Theorem}

\section{Connexion to Prior Literature\label{sec3}}

In this section we discuss the two ingredients developed in this paper:
 Wiener lemma type results, and the faithful tracial state in Gabor analysis.
For each of these we discuss prior results and approaches in literature,
strengths and shortcomings of each method. For precise definitions and
more details of the results we refer the reader to the corresponding
paper.

\subsection{Discussion about Wiener Lemma and Alternate Proofs}

The closest paper to our analysis is \cite{grochJAMS} by Gr\"{o}chenig and 
Leinert. There the authors proved the analogous statement to
Theorem \ref{t1} but only for the algebra of time-frequency shifts
from a lattice. As an off-shot of that result, the authors obtained a
very nice localization result regarding dual Gabor frame generators.
More specifically, if $\{g_{m,n;\al,\be}~;~m,n\in\Z^d\}$ is a 
Gabor frame for $L^2(\R^d)$ with $g\in M^1_v$, the {\em modulation space}
associated to an admissible weight $v$, then the canonical dual frame
$\{\tilde{g}_{m,n;\al,\be}~;~m,n\in\Z^d\}$ is so that $\tilde{g}\in M^1_v$.
One may ask whether the methods used in that paper hold in our case. 
The answer is affirmative. Indeed, the main tool used in \cite{grochJAMS}
is the fact that the Banach algebra $l^1_v(\Z^{2d})$, with the twisted
convolution $\sharp$, is symmetric which in turn is a consequence of 
Hulanicki's and Ludwig's results on symmetry of group algebras associated to 
nilpotent groups. The very same result applies to our case where the
discrete countable group $\Z^{2d}$ is simply replaced by $\R^{2d}$.
Hence, as mentioned in introduction, Theorems \ref{t1} and \ref{t2} 
are consequences of Hulanicki's Theorem from e.g. 
\cite{Hulanicki70}. However, by invoquing a general abstract result one
does not obtain the norm estimates of Theorem \ref{t3} nor the localization
results of Theorems \ref{t5} and \ref{t6}. We prefered to present an
explicit and self-contained proof of Theorem \ref{t1} for two reasons:
(i) for the benefit of reader unfamiliar with the symmetry of nilpotent
group algebras; and (ii) to obtain explicit estimates of the inverse norm.
All these being said we do not want to diminish in any way the beauty
and ellegance of \cite{grochJAMS}.  

We mentioned earlier an important consequence contained in \cite{grochJAMS}, 
namely the localization result of the canonical
dual Gabor frame generator. One may ask if there is an analogous 
consequence in our more general case. The most natural guess would be
to look at irregular Gabor frames and analyze its canonical dual. 
Unfortunately, unlike the regular (i.e. lattice) case, the frame operator
may not necessarily belong to $\Ac_v$ and hence no conclusion can be
drawn from our analysis. Fortunately one can
use another approach to recover the results of \cite{grochJAMS} and 
prove the localization result in the irregular case. This alternative
approach is used in \cite{excess3} and is based on Baskakov's result
mentioned earlier in the introduction. Indeed, the frame operator of
an irregular Gabor frame with generator in $M^1$ 
has a matrix representation with respect to a ``nice'' Gabor frame that
is dominated by a Toeplitz matrix with $l^1$ generating sequence. 
In \cite{excess3} such frames are called 
{\em $l^1$ localized frames}  (more specifically
with respect to a ``nice'' reference Gabor frame generated
by the Gaussian window). The associated matrix of such frames admits a 
pseudoinverse, because of frame condition. Using Baskakov Theorem and 
holomorphic functional calculus one obtains that the pseudoinverse has
the same off-diagonal decaying property which proves the localization
result for the canonical dual frame. In the regular case, the inverse of
an invertible operator that is a linear combination of time-frequency shifts
from a lattice is also a linear combination of time-frequency shifts of the
same lattice. (Here we use ``linear combination'' to denote the generators
of a $C^*$ algebra, hence convergence in operator norm).  
Thus distinct time-frequency labels associated to the inverse operator
are always well separated. This fact combined with Baskakov's result applied
to the pseudoinverse matrix gives an alternative proof to the case considered
in \cite{grochJAMS} (see also \cite{GrochLeinert05}).

The irregular case is fundamentally different from the lattice case, although
it is true that an operator in $\Ac$ has support always contained into a 
countably generated discrete group of the time-frequency plane.
However the main obstrauction in the irregular case is the 
fact that the time-frequency labels of the
inverse of an operator $T\in\Ac$ are not necessarily well-separated, even when
$T$ has finite support. Indeed, if $T=\sum_{\lambda\Lambda}c_\lambda 
U_\lambda$ with $|\Lambda|<\infty$ is invertible in $B(L^2(\R^d))$,
then $T^{-1}=\sum_\mu d_\mu U_\lambda$ with convergence in operator norm 
in $\Ab$. However, in general, $supp(T^{-1})$ has accumulation points in 
$\R^d$. This fact makes difficult the application of Baskakov's Theorem
to irregular frames. To better understand this obstruction, we remark
here only that the conclusion that can be drawn along this line of reasoning
is the following statement. If $T=\sum_{\lambda}c_\lambda U_\lambda\in \Ac$
is invertible in $B(L^2(\R^d))$ then its inverse $T^{-1}=\sum_\mu d_\mu U_\mu$
satisfies
\begin{equation}
\sum_{k\in\Z^{2d}} \sup_{\mu\in B_1(k)}|d_\mu| <\infty
\end{equation}
Clearly this statement is weaker than Theorem \ref{t1} that claims
$\sum_{\mu}|d_\mu| < \infty$.

\subsection{Faithful Tracial States in Gabor Analyis}

In \cite{dalala} Daubechies, Landau, and Landau computed explicitely 
the faithful tracial state on the $W^*$ algebra $\Wb_{a,b}$ generated by 
$\{M_{2\pi ma}T_{nb}~;~m,n\in\Z\}$. They showed that $\Wb_{a,b}$ is a 
$II_1$ factor for $ab\not\in\Q$ (result also known from the rotation
algebra theory, see e.g. \cite{rieffel}), that has a unique faithful
tracial state. In general, for arbitrary $a,b>0$ a faithful
tracial state is defined as the coefficient $c_{0,0}$ of its strongly
convergent, uniquely defined decomposition 
$T=\sum_{m,n}c_{m,n}U_{nb,2\pi ma}$. 
They showed this number ($C_{0,0}$) is computable using the formula
\begin{equation}
\label{eq:trDALALA}
c_{0,0}=\frac{1}{ab}\sum_{k=0}^J \ip{T 1_{I_k}}{1_{I_k}}
\end{equation}
where $J$ is the largest integer smaller than or equal to $ab$, and
the $J+1$ intervals $I_0$, $I_1$, ... ,$I_{J-1}$, $I_J$ are given by
$[0,\frac{1}{a}]$, $[\frac{1}{a},\frac{2}{a}]$, ..., $[\frac{J-1}{a},
\frac{J}{a}]$,$[\frac{J}{a},b]$. 

In this paper we extend this tracial state from algebra generated by
time-frequency shifts from a lattice to the algebra generated by all
time-frequency shifts. Note however the following limitation of our
method. In \cite{dalala} the faithful tracial state applies to a $W^*$
algebra, whereas our Theorem \ref{t3} applies only to a $C^*$
algebra.  The tracial state $\gamma$ of (\ref{eq:3.0}) cannot be
extended to the $W^*$ algebra generated by $\Ac$ since this $W^*$ algebra
is the entire algebra $B(L^2(\R^d))$ (which does not admit a faithful tracial
state). Consider now the series of $C^*$ algebras $(\Cb_{a,b})_{a,b>0}$
each generated by respectively $\{U_{nb,2\pi ma}~;~m,n\in\Z\}$ (we
restrict ourselves here to the one-dimensional case for convenience of 
comparison). For any $T\in \Cb_{a,b}$, for some $a,b>0$, its trace
$\gamma(T)$ can be computed either by (\ref{eq:trDALALA}), or 
by (\ref{eq3.1}). Our formula (\ref{eq3.1}) has the advantage of being
independent of lattice parameters $(a,b)$. In particular this
shows the tracial states defined by (\ref{eq:trDALALA}) are compatible
on operators that belong simultaneously to two different $W^*$ algebras
 (for instance $T\in\Cb_{2a,2b}\subset \Cb_{2a,b}\cap\Cb_{a,2b}$).


We end this section with a comment on Theorem \ref{t4}. In this paper
we solve a restricted case of the HRT conjecture, namely we rule out
the existence of isolated eigenvalues of finite multiplicity for all
finite linear combinations of time-frequency shifts. In fact we obtain
this conclusion for any operator of $\Ac$, hence also for infinite
linear combinations of time-frequency shifts with coefficients in
$l^1$. The other case that was ruled out is the lattice case, that is
when the finitely many time-frequency shifts are from a lattice. This
was beautifully proved by Linnel in \cite{linnel}. One may ask whether
the same arguments hold in our more general case.  There is a
difficulty in trying to do so, namely $\gamma$ is a faithful tracial
state on a $C^*$ algebra (in this paper) unlike the $W^*$ algebra $\Wb_{a,b}$
considered in \cite{linnel}. This difference prevents us from having a
similar proof in our setting. We currently study ways to bypass
this difficulty.

\section{Proof of Results\label{sec4}}

The order of proofs is the following. First we prove the spectral invariance
Theorem \ref{t2}, from where we derive Theorem \ref{t1}.
In the process of proving Theorem \ref{t2} we obtain the norm estimate
(\ref{eq:Tn2}) that allows to prove Theorem \ref{t2.1}.
The corollaries \ref{c0}, \ref{c1} and \ref{c2} follow directly  from \ref{t1}.
The spectral invariance Theorem \ref{t2} can be deduced from the symmetry of
the group algebra of $l^1(\H)$, the Heisenberg group 
$\H=\R^d\times\R^d\times\T$ endowed with a discrete topology, similar
to the approach in \cite{grochJAMS}. In fact those proofs apply almost
verbatim to our case. However we prefer to derive these results directly
in a self-contained manner for two reasons: (i) the proofs are relevant
to researchers not familiar with the group algebra of a nilpotent group,
and (ii) we derive also the explicit norm estimate (\ref{eq:Tn2}) and 
(\ref{eq:Tninv}) which otherwise would not have been available from
the abstract results invoqued in \cite{grochJAMS}.

In Theorem \ref{t3} we construct the faithful tracial state on
$\Ac$, and therefore $\Ab$. This will be proved later in this section.
From this result we will derive Corollaries \ref{c3.1}, \ref{c3},
 and Theorem \ref{t4}. Theorems \ref{t5} and \ref{t6} will follow after
extension to the Banach space $L^{2,\infty}(\R^d\times I)$. 

\subsection{Proof of Theorems \ref{t2}, \ref{t1}, \ref{t2.1} and Corollaries 
\ref{c0},\ref{c1},\ref{c2}}

Theorem \ref{t2} is obtained in two steps. First step involves finite linear
combinations of time-frequency shifts. In the second step we extend the
spectral result to the entire algebra $\Ac$.

Consider $T=\sum_{\lambda\in\Lambda}c_\lambda U_\lambda$ with
$|\Lambda|<\infty$ a finite linear combination of time-frequency shifts.
Note:
\[ T^n= \sum_{\sigma\in\Sigma}d_\sigma U_\sigma \]
where:
\begin{eqnarray}
\Sigma & = & \Lambda + \Lambda + \cdots + \Lambda = \{
\lambda_1+\cdots+\lambda_n~|~\lambda_1,\ldots,\lambda_n\in\Lambda \} \\
d_\sigma & = & \sum_{\begin{array}{c}
\mbox{$\lambda_1,\ldots,\lambda_n\in\Lambda$} \\
\mbox{$\lambda_1+\cdots+\lambda_n=\sigma$} \end{array} }
c_{\lambda_1}\cdots c_{\lambda_n}e^{-it_1\omega_2}e^{-it_2\omega_3}
\cdots e^{-it_{n-1}\omega_{n}}~~,~~{\rm where~each}~
\lambda_k=(t_k,\omega_k)
\end{eqnarray}
Then by Cauchy-Schwarz,
\begin{equation}\label{eq:4.1.2}
\norm{T^n}_{\Ac_v} = \sum_{\sigma\in\Sigma}|d_\sigma|v(\sigma)
\leq |\Sigma|^{1/2}\sup_{\sigma\in\Sigma}|v(\sigma)| \norm{d}_2
\end{equation}
We estimate next the three factors of the right-hand side, and we will
prove the following

\begin{Lemma}\label{l2}
Assume $T=\sum_{\lambda\in \Lambda}c_\lambda U_\lambda$
is an operator in $\Ac$ with $|\Lambda|<\infty$ and 
$\Lambda\subset E_{R_0}(0)$. Then for any submultiplicative weight $v$, that
is $v(x+y)\leq v(x)v(y)$ for all $x,y$, 
\begin{equation}\label{eq:Tn2}
\norm{T^n}_{\Ac_v} \leq \left( n+1 
\right)^{|\Lambda|/2}w(nR_0)\cdot\norm{T^n}_{B(L^2(\R^d))}
\end{equation}
\end{Lemma}
where $w(a)=sup_{\norm{x}=a}v(x)$.
{\bf Proof of Lemma \ref{l2}}

 The cardinal of set $\Sigma$ in (\ref{eq:4.1.2}) 
is upper bounded as follows. Notice that $\lambda_i+\lambda_j=
\lambda_j+\lambda_i$ therefore any permutation of terms in
$\lambda_1+\cdots+\lambda_n$ would produce the same point $\sigma$. Hence:
\begin{equation}\label{eq:4.1.3}
|\Sigma|\leq |\{(k_1,k_2,\ldots,k_{|\Lambda|}~~|~~
k_1+k_2+\ldots +k_{|\Lambda|}=n~,~ k_1,k_2,\ldots,k_{|\Lambda|}\geq 0\}|
=\frac{(n+|\Lambda|)\cdots(n+1)}{|\Lambda|!} \leq (n+1)^{|\Lambda|}
\end{equation}
For the second factor in (\ref{eq:4.1.2}) we need to estimate the radius $R_n$
of a ball $E_R(0)$ in $\R^{2d}$ that includes all $\Sigma$. If $R_0=max_{
\lambda\in\Lambda}\norm{\lambda}$, then for $R_n=nR_0$ we have $\Sigma\subset
E_{R_n}(0)$. Since the weight $v$ is radially non-decreasing,
\begin{equation}\label{eq:4.1.4}
\max_{\sigma\in\Sigma}v(\sigma) \leq w(nR_0)
\end{equation}
The third factor in (\ref{eq:4.1.2}) is a bit more complicated. 
We need to use the following lemma, which is of intrinsec interest:
\begin{Lemma}\label{l1}
For any finite set of time-frequency points 
$\Sigma=\{\sigma_1,\ldots,\sigma_N\}$ there is a function $g\in L^2(\R^d)$
so that $\{U_{\sigma_k}g~;~1\leq k\leq N\}$ is an orthonormal set.
\end{Lemma}
Assume this lemma is proved. Then we apply to our set $\Sigma=\Lambda+\cdots+
\Lambda$ and we obtain, on the one hand
\[ \norm{T^ng}^2=\norm{\sum_{\sigma\in\Sigma}d_\sigma U_\sigma g}^2=
\sum_{\sigma\in\Sigma}|d_\sigma|^2 = \norm{d}_2^2 \]
and on the other hand
\[ \norm{T^n g}^2 \leq \norm{T^n}^2\norm{g}^2=\norm{T^n}^2 \]
Thus we get:
\begin{equation}\label{eq:4.1.5}
\norm{d}_2 \leq \norm{T^n}
\end{equation}
Putting together (\ref{eq:4.1.3},\ref{eq:4.1.4},\ref{eq:4.1.5}) into
(\ref{eq:4.1.2}) we obtain:
\begin{equation}\label{eq:Tn}
\norm{T^n}_{\Ac_v} \leq (n+1)^{|\Lambda|/2}
 w(nR_0)\norm{T^n}_{B(L^2(\R^d))}
\end{equation}
which proves Lemma \ref{l2}. Q.E.D.
\begin{Remark} Inequality (\ref{eq:4.1.5}) follows also idependently from
(\ref{eq:norms}) of Corollary \ref{c3.1}.
\end{Remark}
Taking the $n^{th}$ root and passing to the limit $n\rightarrow\infty$
we obtain:
\begin{equation}
r(T)_{\Ac_v}\leq r(T)_{B(L^2(\R^d))}
\end{equation}

Since the inclusion $\Ac_v\subset
B(L^2(\R^d))$ implies the inverse inclusion of the spectra 
$sp_{B(L^2(\R^d))}(T)\subset sp_{\Ac_v}(T)$, one obtains
 $r_{B(L^2(\R^d))}(T)\leq
r_{\Ac_v}(T)$. This concludes the proof of the spectral radius equation
(\ref{eq3.3.1}) for finite linear combinations of time-frequency shifts.
Before going to the second step, we prove Lemma \ref{l1}.

{\bf Proof of Lemma \ref{l1}}

The statement is equivalent to finding a function $g\in L^2(\R^d)$ so that
$\ip{U_{\sigma}g}{g}=\delta_{\sigma,0}$, for all $\sigma\in
\Delta:=(\Sigma-\Sigma)$.
Let $\Tau$ be the projection of $\Delta$ on the the first factor $\R^d$,
 and $\Omega$ be the projection onto the second factor $\R^d$. Thus
$\Delta \subset \Tau\times\Omega$. Notice both $\Tau$ and $\Omega$ are
finite sets of points of $\R^d$ symmetric about and containing the origin. 
Let $\Omega\setminus\{0\} = \{\omega_1,\ldots,\omega_M\}$ be an enumeration
of $\Omega$, and let $\tau_{min},\tau_{max}>0$ be the radia of two balls
around the origin in $\R^d$ so that $E_{\tau_{min}}(0)\cap\Tau=\{0\}$ and
$\Tau\subset E_{\tau_{max}}(0)$. We set $g$ as follows:
\begin{equation}
\label{eq:g}
g=\sqrt{h_1 * h_2 * \cdots * h_M} / \norm{\sqrt{h_1 * h_2 * \cdots * h_M}}
\end{equation}
where $*$ denotes the usual convolution, and $g_1,g_2,\ldots,g_M$ are 
constructed as follows. First we construct inductively the sequence 
$t_1,t_2,\ldots,t_M\in\R^d$ so that:
\begin{enumerate}
\item $\ip{t_1}{\omega_1} = (2n_1+1)\pi$ for some integer $n_1\in\Z$ and
$\norm{t_1}>\tau_{max}$
\item Assume $t_1,t_2,\ldots,t_k$ were set; then $t_{k+1}$ is constructed so
that: 
(i) $\ip{t_{k+1}}{\omega_{k+1}}=(2n_{k+1}+1)\pi$ for some integer 
$n_{k+1}\in\Z$, and 
(ii) $\norm{t_{k+1}}>\norm{t_1}+\cdots+\norm{t_k}+2M\tau_{max}$
\end{enumerate}
With this choice for $\{t_1,\ldots,t_M\}$, we set:
\begin{equation}
h_k = 1_{E} + 1_{t_k+E}
\end{equation}
where $E=E_{\tau_{min}/M}(0)$ is the Euclidian ball of radius $\tau_{min}$
centered at the origin, and $1_E$, respectively, $1_{t_k+E}$, is the
characteristic function of $E$, respectively of $t_k+E$. Note that $g$
is a sum of $2^M$ ``bump'' functions each supported inside balls of radius
$\tau_{min}$ and each at a distance from one another of at least $\tau_{max}$.
Thus all translates with shifts from $\Tau\setminus\{0\}$ are disjoint.
 Hence $\ip{U_\mu g}{g}=0$ for all $\mu=(t,\omega)\in\Delta$
 with $t\in\Tau\setminus\{0\}$. It remains to check only that 
$\ip{M_{\omega_k}g}{g}=0$. Using Fourier transform, this is equivalent to
 $\Fc(|g|^2)(\omega_k)=0$. But the choice of $t_k$ guarantees that 
$\Fc(h_k)(\omega_k)=0$ which concludes the proof of Lemma \ref{l1}. Q.E.D.

Now we are ready to go to step 2 of the proof of Theorem \ref{t2}.
Consider now $T=\sum_{\lambda}c_\lambda U_\lambda\in Ac_v$. Fix
$\eps>0$.  Let $\Lambda$ be the finite set so that
$\sum_{\lambda\in\R^{2d}\setminus\Lambda}|c_\lambda|v(\lambda)<\eps$. Set
$T_0=\sum_{\lambda\in\Lambda}c_\lambda U_\lambda$ and $R=T-T_0$. Thus
$\norm{R}\leq\norm{R}_{\Ac_v}<\eps$. For a more convenient notation in
the following we denote by $s=(s_\sigma)_\sigma$ and
$r=(r_\rho)_\rho$ the coefficients of $T_0$, respectively $R$,
\begin{equation}
s_\sigma = \left\{ \begin{array}{rcl}
\mbox{$c_\sigma$} & if & \mbox{$\sigma\in\Lambda$} \\
0 & & otherwise
\end{array} \right. ~~~~
r_\rho = \left\{ \begin{array}{rcl}
0 & if & \mbox{$\rho\in\Lambda$} \\
\mbox{$c_\rho$} & & otherwise
\end{array} \right.
\end{equation}
Next expand $T^n=(T_0+R)^n$ as follows
\begin{eqnarray}
T^n & = & \sum_{m=0}^n\sum_{k,j}R^{k_1}T_0^{j_1}R^{k_2}T_0^{j_2}\cdots
R^{k_l}T_0^{j_l} =\sum_{m=0}^n\sum_{k,j}\sum_{\lambda\in\R^{2d}}\left(
r^{k_1}\sharp s^{j_1}\sharp \cdots \sharp r^{k_L}\sharp s^{j_L}\right)_\lambda
U_\lambda
\end{eqnarray}
where $k=(k_1,k_2,\ldots,k_L)$, $j=(j_1,j_2,\ldots,j_L)$ are vectors of
nonnegative integers so that $k_1+\cdots k_L=n-m$ and $j_1+\cdots+j_L=m$,
 and $r^k=r\sharp\cdots\sharp r$ is the $k$-fold twisted convolution.
Then the $\lambda$-coefficient expands into

\begin{eqnarray}
c_\lambda & = & \sum_{\rho_1,\cdots,\rho_L\in\R^{2d}}
\sum_{\sigma_1,\cdots,\sigma_L\in\R^{2d}}
r^{k_1}_{\rho_1}e^{-i \ip{a_1}{d_1-b_1}} s^{j_1}_{\sigma_1-\rho_1} 
e^{-i \ip{c_1}{b_2-d_1}} r^{k_2}_{\rho_2-\sigma_1}e^{-i\ip{a_2}{d_2-b_2}}
s^{j_2}_{\sigma_2-\rho_2}\cdots \\
& &\cdot e^{-i\ip{c_{L-1}}{b_L-d_{L-1}}} r^{k_L}_{\rho_L-\sigma_{L-1}}
e^{-i\ip{a_L}{d_L-b_L}} s^{j_L}_{\sigma_L-\rho_L} \delta_{\sigma_L,\lambda}
\end{eqnarray}
where $\rho_l=(a_l,b_l)$ and $\sigma_l=(c_l,d_l)$ are the components
 of the $2L$ phase-space points constrained by $\sigma_L=\lambda$ 
as expressed by the last Kronecker term.

The next step is to change the summation variables and rearrange the terms 
as suggested by Hulanicki in \cite{Hulanicki66}. 
Let $\tilde{\rho}_p=\rho_p-\sigma_{p-1}=(\tilde{a}_p,\tilde{b}_p)$, 
$1\leq p\leq L$, with convention $\sigma_0=(0,0)$. Also denote
by $V_p$ the unitary 
$$(V_p s)_{\lambda=(t,\omega)}=
e^{-i\ip{\tilde{a}_p}{\omega}}s_{\lambda-\tilde{\rho}_p}.$$
Then $c_\lambda$ turns into:
\begin{eqnarray}
c_\lambda & = & \sum_{\tilde{\rho}_1,\cdots,\tilde{\rho}_L\in\R^{2d}}
  e^{i a_1b_1}
r^{k_1}_{\tilde{\rho}_1} r^{k_2}_{\tilde{\rho}_2}\cdots 
r^{k_L}_{\tilde{\rho}_L} \sum_{\sigma_1,\cdots,\sigma_L}
(V_1 s^{j_1})_{\sigma_1}e^{-i\ip{c_1}{d_2-d_1}}
(V_2s^{j_2})_{\sigma_2-\sigma_1}e^{-i\ip{c_2}{d_3-d_2}} \\
& & \cdot\ldots  \cdot e^{-i\ip{c_{L-1}}{d_L-d_{L-1}}}(V_L s^{j_{L}})_{
\sigma_L-\sigma_{L-1}}\delta_{\sigma_L,\lambda} \\
& = & \sum_{\tilde{\rho}_1,\cdots,\tilde{\rho}_L\in\R^{2d}}
  e^{i a_1b_1}
r^{k_1}_{\tilde{\rho}_1} r^{k_2}_{\tilde{\rho}_2}\cdots 
r^{k_L}_{\tilde{\rho}_L} ((V_1 s^{j_1})\sharp(V_2 s^{j_2})
\sharp\cdots\sharp(V_L s^{j_L}))_\lambda
\end{eqnarray}
and thus
\begin{equation}
|c_\lambda|\leq \sum_{\tilde{\rho}_1,\ldots,\tilde{\rho}_L}
|r^{k_1}_{\tilde{\rho}_1}|\cdot\ldots\cdot|r^{k_L}_{\tilde{\rho_L}}|\cdot
|(V_1 s^{j_1}\sharp V_2 s^{j_2}\sharp\cdots \sharp V_L s^{j_L})_{\lambda}|
\end{equation}
Since
$$w(\lambda)\leq
w(\tilde{\rho}_1)\cdot\ldots w(\tilde{\rho}_L)\cdot w(\lambda-\tilde{\rho}_1
-\cdots-\tilde{\rho}_L)$$
and 
$$ supp(V_1 s^{j_1}\sharp V_2 s^{j_2}\sharp\cdots\sharp V_L s^{j_L})\subset
(\tilde{\rho}_1+\cdots+\tilde{\rho}_L)+\Sigma~~,~~
\Sigma:=\underbrace{\Lambda+\Lambda+\cdots+\Lambda}_m $$
we obtain
\begin{equation}
\norm{T^n}_{\Ac_v}\leq \sum_{m=0}^n \left( \begin{array}{c}
n\\ m \end{array} \right) \eps^{n-m}  
\max_{\sigma\in\Sigma}v(\sigma)
\sup_{\tilde{\rho}_1,\ldots,\tilde{\rho}_L}
\, \sum_{\lambda}  
\left|((V_1 s^{j_1})\sharp\cdots\sharp (V_L s^{j_L}))_\lambda\right|
\end{equation}
Now we will estimate the sum over $\lambda$
above similar to the estimation in (\ref{eq:4.1.2}). 
The cardinal of $\Sigma$ has been shown in Lemma \ref{l1} to be bounded by 
$(m+1)^{|\Lambda|}$ and hence by $(n+1)^{|\Lambda|}$. 
Note also that $\Sigma\subset E_{mR_0}(0)$ where $R_0$
is a radius so that $\Lambda\subset E_{R_0}(0)$. Thus we get:
\begin{equation}\label{eq:4.44}
\norm{T^n}_{\Ac_v} \leq (n+1)^{|\Lambda|/2} w(nR_0)
\sum_{m=0}^n \left( \begin{array}{c}
n\\ m \end{array} \right) \eps^{n-m} 
\sup_{a_1,\ldots,a_L,b_1,\ldots,b_L} \left( \sum_{\lambda}
|((V_1 s^{j_1})\sharp\cdots\sharp (V_L s^{j_L}))_\lambda|^2 \right)^{1/2}
\end{equation}
By Lemma \ref{l2} the $l^2$ norm of the sequence $s=(V_1 s^{j_1}\sharp\cdots
\sharp V_L s^{j_L})$ is bounded by the operator norm obtained by linear 
combinations of time-frequency shifts with coefficients from $s$:
\[ \norm{s}_2 \leq \norm{\sum_{\lambda}s_\lambda U_\lambda}_{B(L^2(\R^d))} \]
Note the operator associated to $s$ is (up to a constant phase factor):
\[ U_{\tilde{\rho}_1}T_0^{j_1} U_{\tilde{\rho}_2}T_0^{j_2}\cdot\ldots\cdot
U_{\tilde{\rho}_L}T_0^{j_L} \]
Thus we get:
\begin{equation}
\norm{s}_2 \leq \norm{T_0}^{j_1+\cdots+j_L}_{B(L^2(\R^d))}=\norm{T_0}^m_{B(
L^2(\R^d))}
\end{equation}
which turns (\ref{eq:4.44}) into:
\begin{equation}
\norm{T^n}_{\Ac_v}\leq (n+1)^{|\Lambda|}(\eps+\norm{T_0}_{B(L^2(\R^d))})^n
\end{equation}
Now taking the $n^{th}$ root and passing to the limit $n\rightarrow\infty$
we obtain:
\[ r_{\Ac_v}(T)\leq \eps+\norm{T_0}_{B(L^2(\R^d))} \]
Since $\eps>0$ was arbitrary, and $\norm{T_0}\leq\norm{T}$ we obtain:
\begin{equation}
\label{eq:4.45}
r_{\Ac_v}(T)\leq \norm{T}_{B(L^2(\R^d))}
\end{equation}
Since $r_{\Ac_v}(T^n)= (r_{\Ac_v}(T))^n$ we obtain:
\[ r_{\Ac_v}(T) \leq \norm{T^n}^{1/n}_{B(L^2(\R^d))} \]
and passing to the limit $n\rightarrow\infty$ we obtain
\begin{equation}
\label{eq:4.46}
r_{\Ac_v}(T)\leq r_{B(L^2(\R^d))}(T)
\end{equation}
The converse inequality is immediate from $\Ac_v\subset B(L^2(\R^d))$. 
This ends the proof of (\ref{eq3.3.1}). 
Theorem \ref{t2}. $\Diamond$
\vspace{5mm}

Now Theorem \ref{t1} is immediate. For completeness we include its proof.

{\bf Proof of Theorem \ref{t1}}

Assume $T\in\Ac_v$ is invertible in $B(L^2(\R^d))$. Then there are finite
 $A,B>0$ so that $A\leq T^*T\leq B$. Thus we have:
\[ T^{-1}=(T^*T)^{-1}T^* = \frac{2}{A+B}(1-(1-\frac{2}{A+B}T^*T))^{-1}T^* \]
Let $R=1-\frac{2}{A+B}T^*T$. Note $\norm{R}_{B(L^2(\R^d))}
=\frac{B-A}{B+A}<1$, hence 
\begin{equation}
\label{eq:80}
 T^{-1} = \frac{2}{A+B}\sum_{n\geq 0}(1-\frac{2}{A+B}T^*T)^nT^*
\end{equation}
But by Theorem \ref{t2}, 
\[ \lim_{n\rightarrow\infty}
(\norm{R^n}_{\Ac_v})^{1/n}=r_{\Ac_v}(R)=r_{B(L^2(\R^d))}(R)=
\norm{R}_{B(L^2(\R^d))}=\frac{B-A}{B+A}<1. \]
which proves the series $\sum_{n\geq 0}R^n$ is convergent in $\Ac_v$ hence
$T^{-1}\in\Ac_v$. Q.E.D.

During the proof of Theorem \ref{t2} we obtained the
estimate (\ref{eq:Tn2}). We will use this in proving Theorem \ref{t2.1}.

{\bf Proof of Theorem \ref{t2.1}}

Assume $T$ is invertible in $B(L^2(\R^d))$. Then for
 $A=\norm{T^{-1}}^{-2}$ and $B=\norm{T}^2$,
\[ 0 < A\leq T^*T \leq B <\infty \]
Note $\norm{1-\frac{2}{A+B}} = \frac{B-A}{B+A}<1 $. Thus
\[ T^{-1} = (T^*T)^{-1}T^*=\frac{2}{A+B}\sum_{n\geq 0}
(1-\frac{2}{A+B}T^*T)^n T^* \]
which converges in operator norm in $B(L^2(\R^d))$.
The estimate (\ref{eq:Tn2}) of Lemma \ref{l2} turns into:
\[ \norm{1-\frac{2}{A+B}T^*T}_{\Ac_v}\leq \left(n+1
 \right)^{|\Lambda'|/2} w(nR_0') \norm{1-\frac{2}{A+B}T^*T}^n \]
where $\Lambda'$ is the label set of $1-\frac{2}{A+B}T^*T$, and $R_0'$ is
 so that $\Lambda'\subset E_{R_0'}(0)$. Since $\Lambda'\subset\Lambda-\Lambda$
we have $|\Lambda'|\leq 2|\Lambda|=2N$ and $R_0'\leq 2R_0$, 
where $R_0=\max_{\lambda\in\Lambda}\norm{\lambda}$ is so that
$\Lambda\subset E_{R_0}(0)$. Thus we get
\[ \norm{T^{-1}}_{\Ac_v}\leq \frac{2}{A+B}
\sum_{n\geq 0}\left( n+1 \right)^{N}
w(2nR_0)\left(\frac{B-A}{B+A}\right)^n \norm{T}_{\Ac_v}\]
For $w(x)=C(1+x)^m$ and $\rho=max(1,2R_0)$, $w(2nR_0)\leq \rho^m(1+n)^m$
and for $\theta_0=(B-A)/(B+A)$ we obtain
\[ \norm{T^{-1}}_{\Ac_v}\leq \frac{2C\rho^m\norm{T}_{\Ac_v}}{A+B}\sum_{n\geq 0}
(1+n)^{m+N}\theta^n \leq \frac{2C\rho^m\norm{T}_{\Ac_v}}{A+B} 
[\frac{d^{m+N}}{d\theta^{m+N}}\sum_{n\geq 0}\theta^n]{|}_{\theta=\theta_0} \]
Since $\theta_0<1$, by direct summation and then differentiation
we obtain (\ref{eq:Tninv}) which
ends the proof of this Theorem. Q.E.D.

\subsection{Proof of Theorems \ref{t3}, \ref{t4} and Corollaries \ref{c3.1} and
\ref{c3}}

For any bounded operator $T\in B(L^2(\R^d))$ denote by $a_{M,N}(T)$ the
following expression
\begin{equation}
a_{M,N}(T) = \frac{1}{(\al\be)^d(2M+1)^d(2N+1)^d}\sum_{|m|\leq M}
\sum_{|n|\leq N}\ip{Tg_{m,n;\al,\be}}{\tilde{g}_{m,n;\al,\be}}
\end{equation}

{\bf Trace on $\Ac$}

First we need to show that for every $T\in\Ac$
the limit $\lim_{M,N\rightarrow\infty}a_{M,N}$ exists and equals 
$c_0=\gamma(T)$, the $0$-coefficient of $T$. 
We prove this statement in two steps. First we consider the unitary generators
of $\Ac$, and then we extend by continuity to the entire $\Ac$.

\begin{Lemma}\label{l4.1}
Let $U_\lambda$ denote the time-frequency shift with parameter 
$\lambda=(t,\omega)$. Then
\begin{equation}
\label{eq:traceUlambda}
\lim_{M,N\rightarrow\infty}a_{M,N}(U_\lambda)=\left\{
\begin{array}{rcl}
1 & if & \mbox{$\lambda=0$} \\
0 & if & \mbox{$\lambda\neq 0$}
\end{array} \right.
\end{equation}
\end{Lemma}

{\bf Proof}

We explicitely compute $a_{M,N}(U_\lambda)$
\begin{equation}\label{eq:pert}
 a_{M,N}(U_\lambda) = \frac{1}{(\al\be)^2(2M+1)^d(2N+1)^d}
\sum_{|m|\leq M}\sum_{|n|\leq N} e^{2\pi im\al t}e^{in\be \omega}
\int e^{i\omega x} g(x-t)\overline{\tilde{g}(x)} dx
\end{equation}
There are now two cases:

Case 1. $(t,\omega)=(\frac{K}{\al},\frac{2\pi J}{\be})$ for some $K,J\in\Z^d$.
Then $e^{2\pi im\al t}=e^{2\pi in\be}=1$ and summations over $m$ and $n$
cancel $(2M+1)^d(2N+1)^d$ factor
\[ a_{M,N}(U_\lambda) = \int e^{\frac{2\pi i}{\be}J x}g(x-K\al)
\overline{\tilde{g}(x)}\,dx \]
Recall $\Gc$ is a Gabor frame with dual frame $\tilde{\Gc}$. By duality
principle (\cite{dalala},\cite{Janss},\cite{RonShen}) 
$\Gc'=\{g_{m,n;\frac{1}{\be},\frac{1}{\al}}~;~m,n\in\Z^d\}$ and
 $\tilde{\Gc}'=\{\frac{1}{(\al\be)^d}\tilde{g}_{m,n;\frac{1}{\be},
\frac{1}{\al}}~;~m,n\in\Z^d\}$
are Riesz basic sequences biorthonormal to one another. Thus
\[ \ip{g_{J,K;\frac{1}{\be},\frac{1}{\al}}}{\tilde{g}} = (\al\be)^d\delta_{J,0}
\delta_{K,0} \]
and combined with (\ref{eq:pert}) proves (\ref{eq:traceUlambda}) in this
case.

Case 2. $(t\al,\frac{\omega\be}{2\pi})\not\in\Z^{2d}$. Then for
 $t\al\not\in\Z^d$ a direct computation shows
\[ \lim_{M\rightarrow\infty}\frac{1}{(2M+1)^d}\sum_{|m|\leq M}
e^{2\pi i m\al t} = 0 \]
whereas for $\frac{\omega\be}{2\pi}\not\in\Z^d$
\[ \lim_{N\rightarrow\infty}\frac{1}{(2N+1)^d}\sum_{|n|\leq N}
e^{in\be \omega}=0 \]
This ends the proof of Lemma \ref{l4.1}. Q.E.D.

By linearity we extend the result of this lemma to finite linear combinations
of unitaries $U_\lambda$'s:
\[ \lim_{M,N\rightarrow\infty} a_{M,N}(\sum_{k=1}^L c_{\lambda_k}U_{\lambda_k})
 = \left\{ \begin{array}{rcl}
\mbox{$c_0$} & if & \mbox{$0\in\{\lambda_1,\ldots,\lambda_L\}$} \\
0 & & otherwise 
\end{array} \right. \]
Next the limit extends to the entire $\Ac$ by Lebesgue's domainated
convergence theorem:
\begin{equation}
\lim_{M,N\rightarrow\infty} a_{M,N}(\sum_{\lambda}c_\lambda U_\lambda)
=\sum_{\lambda}c_\lambda \lim_{M,N\rightarrow\infty}a_{M,N}(U_\lambda) = c_0
=\gamma(\sum_\lambda c_\lambda U_\lambda)
\end{equation}
Consider $T=\sum_\lambda c_\lambda U_\lambda\in\Ac$ and $S=\sum_\mu d_\mu U_\mu
\in\Ac$. Then $TS = \sum_\rho (c\sharp d)_\rho U_\rho\in\Ac$ and 
$ST=\sum_\rho (d\sharp c)_\rho U_\rho\in\Ac$. But now
\begin{equation}
\gamma(TS)=(c\sharp d)_0 = \sum_{\lambda}c_\lambda d_{-\lambda} =
(d\sharp c)_0 = \gamma(ST)
\end{equation}
This shows $\gamma$ is a tracial state on $\Ac$.
Finally, for $T=\sum_{\lambda}c_\lambda U_\lambda\in\Ac$. Then
\begin{equation}\label{eq:4.4.4}
 \gamma(T^*T) = \sum_{\lambda}|c_\lambda|^2\geq 0~~,~~\gamma(T^*T)=0~
{\rm iff}~T=0 
\end{equation}
Thus $\gamma$ is a faithful state.

{\bf Extension to $\Ab$}

For $T\in\Ac$ from (\ref{eq3.1}) we obtain:
\[ |\gamma(T)| \leq \norm{T}_{B(L^2(\R^d))} \]
Then we can extend $\gamma$ to the completion of $\Ac$ with respect to the
operator norm. The completion of $\Ac$ is denoted $\Ab$ and is a $C^*$ 
algebra. On this algebra, $\gamma$ remains a faithful tracial state.
 This ends the proof of Theorem \ref{t3}. Q.E.D.

{\bf Corollaries of Theorem \ref{t3}}

Corollary \ref{c3.1} proof is immediate. In particular
 equation (\ref{eq:norms}) follows as in (\ref{eq:4.4.4}),
 whereas (\ref{eq:coeffs}) is a consequence of Lemma \ref{l4.1}.

Corollary \ref{c3} follows as follows. First note that if $T$ is a 
finite rank operator then $trace(T^*T)<\infty$. But then
\[ trace(T^*T) = \sum_{m,n\in\Z^d} \ip{T^*
Tg_{m,n;\al,\be}}{\tilde{g}_{m,n;\al,\be}} < \infty
\]
and thus $\lim_{M,N\rightarrow\infty}a_{M,N}(T^*T) = 0$.
By Corollary \ref{c3.1} this implies $T=0$. Now consider $T$ a compact operator
in $\Ab$. Then $T^*T\in\Ab$ is a non-negative compact operator. Let $s>0$ be
an eigenvalue of $T^*T$ and let $P_s$ denote the projection onto its 
eigenspace. On the one hand $P_s$ is finite rank, since $T^*T$ is compact,
on the other hand, by holomorphic functional calculus (see \cite{RieszNagy}), 
$P_s$ belongs to $C^*$ algebra $\Ab$. Then, as shown before, $P_s$ has to 
vanish, which proves $T=0$. Q.E.D.

{\bf Proof of Theorem \ref{t4}}

Theorem \ref{t4} is a consequence of the holomorphic functional calculus. 
Assume $T=\sum_{\lambda}c_\lambda U_\lambda$ has a finite isolated eigenvalue
say $\mu_0$. Since it is isolated, by holomorphic functional calculus 
(see \cite{RieszNagy}) the
orthogonal projection onto the eigenspace is given by
\begin{equation}
\label{eq:proj}
P_{\mu_0}=\frac{1}{2\pi i}\int_{\Gamma}(zI-T)^{-1}dz
\end{equation}
where $\Gamma$ is a circle in complex plan centered at $\mu_0$ so that it
separates $\mu_0$ from the rest of the spectrum of $T$. Thus 
$P_{\mu_0}\in\Ab$. Since $\mu_0$ has finite multiplicity it follows that
$P_{\mu_0}$ has finite rank but then by Corollary \ref{c3}, $P_{\mu_0}=0$
which ends the proof of Theorem \ref{t4}. Q.E.D.

\subsection{Proof of Theorems \ref{t5} and \ref{t6}}

Theorems \ref{t5} and \ref{t6} characterize finite and half-compactly 
supported operators in $\Ac$. The proof of Theorem \ref{t5} is based on
a spectral radius computation done in Lemma \ref{l2}. In turn, Theorem
\ref{t5} allows the operator extension to the Banach space $L^{2,\infty}$
introduced in Section \ref{sec2b}. Once this extension is established,
Theorem \ref{t6} follows easily.

For a $\rho>0$ we set $f(\lambda)=e^{\rho\norm{\lambda}}$, $\lambda\in\R^d$.
For convenience we denote $\Bc_\rho=\Ac_f$, the Banach algebra of bounded
operators in $\Ac$ whose coefficients decay exponentially fast with rate 
$\rho$.
Note the spectral radius of an operator $T\in\Bc_\rho$ is not the same as the
spectral radius in $B(L^2(\R^d))$, since $f$ does not satisfy the GRS 
condition.

{\bf Proof of Theorem \ref{t5}}

Assume $T=\sum_{\lambda\in\Lambda}c_\lambda U_\lambda$ is an invertible
operator on $L^2(\R^d)$ with finite support. Let $A,B>0$ be the bounds in
\[ A\norm{f} \leq \norm{Tf} \leq B\norm{f} ~~,~~\forall f\in L^2(\R^d) \]
Then , as shown in (\ref{eq:80})
\begin{equation}\label{eq:81}
 T^{-1} = \frac{2}{A+B} \sum_{n\geq 0}R^n T^* 
\end{equation}
where $R=1-\frac{2}{A+B}T^*T$. Since $R$ has finite support let
$N:=|supp(R)|<\infty$ and $R_0>0$ so that $supp(R)\subset E_{R_0}(0)$.

The goal is to show there exists a $\rho>0$ so that $T^{-1}\in\Bc_\rho$.
Clearly each term in (\ref{eq:81}) belongs to $\Bc_\rho$. The only problem
is to check the series converges in $\Bc_\rho$. To do so it is sufficient
to show that
\begin{equation}
\label{eq:82}
\limsup_{n\rightarrow\infty} (\norm{R^n}_{\Bc_\rho})^{1/n} < 1
\end{equation}
Let us choose $\rho = \frac{1}{R_0}ln\frac{B+A}{2(B-A)}$. Then
 $e^{\rho R_0}\norm{R} \leq \frac{1}{2}$. Now we apply Lemma \ref{l2}
to the algebra $\Bc_\rho$, since $f(x)=e^{\rho\norm{x}}$ is a submultiplicative
weight. We obtain:
\[ (\norm{R^n}_{\Bc_\rho})^{1/n} \leq (n+1)^{N/(2n)}e^{\rho R_0} \norm{R} \]
By choice of $\rho$ we obtain (\ref{eq:82}) which ends the proof of 
Theorem \ref{t5}. Q.E.D.

Let us turn now to Theorem \ref{t6}. First we recall the mixed norm 
Banach space introduced in Section. \ref{sec2b}.  For a
compact neighborhood of the origin $I$ we denote by $L^{2,\infty}
(\R^d\times I)$ the mixed norm Banach space defined in (\ref{eq:L2inf}).
An operator $T=\sum_{\lambda\in\Lambda} c_\lambda U_\lambda\in\Ac$ 
extends from $L^2(\R^d)$ to $L^{2,\infty}(\R^d\times I)$ by
\begin{equation}
\label{eq:Text}
f\in L^{2,\infty}(\R^d\times I)~\mapsto ~\Tb f(x,y) = 
\sum_{\lambda=(t,\omega)\in\Lambda}
 c_\lambda e^{i\ip{\omega}{x+iy}}f(x-t,y)~~,~~\forall (x,y)\in\R^d\times I
\end{equation}
Assume that for all $\lambda=(t,\omega)\in\Lambda$, 
$\norm{\omega}\leq \Omega$, in other words $\Lambda\subset\R^d\times 
E_{\Omega}(0)$. Assume also $I\subset E_{\rho}(0)$, for some $\rho>0$.
Then
\[ \norm{\Tb f}_{2,\infty}\leq \sum_{\lambda\in\Lambda}
sup_{y\in I} e^{-\ip{\omega}{y}}(\int_{\R^d}|f(x,y)|^2dx)^{1/2}
\leq e^{\rho\Omega}\sum_{\lambda}|c_\lambda|\norm{f}_{2,\infty} \]
Thus
\begin{equation}\label{eq:Trho}
\norm{\Tb}_{B(L^{2,\infty}(\R^d\times I))}\leq e^{\rho\Omega} \norm{T}_{\Ac}
\end{equation}
This proves $T$ extends to a bounded operator $\Tb$ on 
$L^{2,\infty}(\R^d\times I)$.
Let $T^{(y)}$ denotes the restriction of $\Tb$ to the ``slice'' indexed by $y$
of $L^{2,\infty}$, that is:
\begin{equation}\label{eq:6.1}
T^{(y)}:L^2(\R^d)\rightarrow L^2(\R^d)~~,~~T^{(y)}f(x)=
\sum_{\lambda=(t,\omega)} c_\lambda e^{-\ip{\omega}{y}} U_\lambda f(x)
\end{equation}
Note for $\Lambda\in \R^d\times E_{\Omega}(0)$ each $T^{(y)}\in\Ac$ and
\begin{equation}\label{eq:6.2}
\norm{\Tb}_{L^{2,\infty}(\R^d\times I)} = 
\sup_{y\in I}\norm{T^{(y)}}_{B(L^2(\R^d))}
\end{equation}
Now we are ready to prove Theorem \ref{t6}.

{\bf Proof of Theorem \ref{t6}}

(1) The estimate (\ref{eq:3.6}) follows from (\ref{eq:Trho}).

(2) From (\ref{eq:6.1}), (\ref{eq:6.2}) and Corollary \ref{c3.1} we get:
\[ e^{-\ip{\omega}{y}}|c_\lambda| \leq \norm{T^{(y)}}_{B(L^2(\R^d))}
\leq \norm{T}_{B(X_\rho)} \]
Combine now with (\ref{eq:3.6}) to obtain:
\[ |c_\lambda| \leq C e^{\rho\Omega + \ip{\omega}{y}} \]
Assume there is a $\lambda=(t,\omega)\in\Lambda$ so that 
$\norm{\omega}>\Omega$. Choose $y\in E_\rho(0)$ so that $\ip{\omega}{y}
=-\rho \norm{\omega}$. Then
\[ |c_\lambda| \leq C e^{-\rho(\norm{\omega}-\Omega)} \]
should hold valid for all $\rho>0$. At limit $\rho\rightarrow\infty$ we
get $|c_\lambda|\leq 0$ which proves that 
$supp(T)\subset\R^d\times E_\Omega(0)$.

(3) and (4) follows from (1) and (2) by noticing that Fourier transform
is an intertwining operator between the translation shift and the 
modulation shift that switches the modulation and translation parameters.
 Thus (3) and (4) reduce to (1) and (2). Q.E.D.

\section{Acknowledgments}

The author wishes to thank Thomas Strohmer for pointing out, firstly
the Sj\"{o}strand's paper \cite{sjostrand95}, secondly the Gohberg's
paper \cite{gohberg89}, and finally for comments on using the approach
here to a class of pseudodifferential operators. The trace results
would not have been possible without the stimulating discussions the
author had with Zeph Landau. He is very happy to acknowledge and thank
for all the comments and suggestions received during the past seven
years from him. Many thanks go to Chris Heil for discussions and ideas
exchanged during all these years. Special thanks go to Charly
Gr\"{o}chenig for discussions on the approach used in
\cite{grochJAMS}, and to Hans Feichtinger for comments on the
draft. The author also thanks Hans Feichtinger, Charly Gr\"{o}chenig,
John Benedetto, and Erwin Schr\"{o}dinger Institute for the invitation
to and hospitality during his stay in Vienna at ESI in May 2005, when
most of the proofs of this paper were obtained.

\bibliographystyle{alpha}
\bibliography{analysis}

\end{document}